\documentclass[epsfig,amstex,amssymb]{amsart}
\usepackage{epsfig}

\begin{document}

\title[The Chord Problem]{The Chord Problem and a new method of filling by pseudoholomorphic curves}
\author{
        Casim Abbas     
}
\address{Casim Abbas\\Department of Mathematics\\Michigan State University\\Wells Hall\\East Lansing, MI 48824\\USA}

\date{\today}        
\footnote{This material is based upon work supported by the National Science Foundation under Grant No. 0196122 and by a New York University Research Challenge Fund Grant. This paper was written while the author visited ETH Z\"{u}rich. He would like to thank FIM and Prof. Eduard Zehnder for their hospitality.}  
\begin{abstract}
Let $M$ be a closed three dimensional manifold with contact form $\lambda$ so that $\ker\lambda$ is tight. In this paper we will present a first application of the filling method by pseudoholomorphic curves recently developed by the author. We will show that Legendrian knots ${\mathcal L}\subset M$ satisfying suitable assumptions admit a Reeb Chord, i.e. there is a trajectory $x(t)$ of the Reeb vector field and $T>0$ such that $x(0),x(T)\in{\mathcal L}$ and $x(0)\neq x(T)$.
\end{abstract}
\maketitle

%
%
\parindent0ex
\parskip1ex plus0.4ex minus0.2ex

\newtheorem{lemma}{Lemma}[section]
\newtheorem{theorem}[lemma]{Theorem}
\newtheorem{corollary}[lemma]{Corollary}
\newtheorem{proposition}[lemma]{Proposition}
\newtheorem{definition}[lemma]{Definition}
\newtheorem{conjecture}{Conjecture}
\newtheorem{example}[lemma]{Example}
\newtheorem{condition}[lemma]{Condition}
\newcommand{\degree}{\mbox{deg}}
\newcommand{\ind}{\mbox{ind}}
\newcommand{\abs}{|}
\newcommand{\ve}{\varepsilon}
\newcommand{\Id}{\mbox{Id}}
\newcommand{\GL}{\mbox{GL}}
\newcommand{\cal}{\mathcal}
\newcommand{\eps}{\varepsilon}
\newcommand{\To}{\longrightarrow}
\newcommand{\Real}{{\bf{R}}}
\newcommand{\Complex}{{\bf{C}}}
\newcommand{\RM}{{\bf{R}}\times M}
\newcommand{\RL}{{\bf{R}}\times {\mathcal L}}
\newcommand{\tu}{\tilde{u}}
\newcommand{\pil}{\pi_{\lambda}}
\newcommand{\pas}{\partial_s}
\newcommand{\pat}{\partial_t}
\newcommand{\od}{\{0\}\times{\mathcal D}^{\ast}}
\newcommand{\SC}{{\mathcal S}}
\newcommand{\ho}{\mbox{Hom}}
\newcommand{\oo}{{\mathcal O}}

\section{Introduction}
Let $M$ be a closed $(2n+1)$--dimensional manifold with contact form $\lambda$, i.e. $\lambda$ is a 1--form on $M$ such that $\lambda\wedge (d\lambda)^n$ is a volume form. The contact structure associated to $\lambda$ is the $(2n)$--dimensional vector bundle $\xi=\ker\lambda\rightarrow M$, which is a symplectic vector bundle  with symplectic structure $d\lambda|_{\xi\oplus\xi}$. There is a distinguished vector field associated to a contact form, the Reeb vector field $X_{\lambda}$, which is defined by the equations 
\[
i_{X_{\lambda}}d\lambda\equiv 0\ ,\ i_{X_{\lambda}}\lambda\equiv 1.
\]
The main result of this paper is about the global dynamics of the Reeb vector field on three dimensional contact manifolds. More precisely, we will prove an existence result for so--called 'characteristic chords'. These are trajectories $x$ of the Reeb vector field which hit a given Legendrian submanifold ${\mathcal L}$ at two different times $t=0$ and $T>0$. We also ask for $x(0)\neq x(T)$, otherwise the chord would actually be a periodic orbit. Recall that an $n$--dimensional submanifold ${\mathcal L}$ of a $(2n+1)$--dimensional contact manifold $(M,\xi)$ is called Legendrian if it is everywhere tangent to the hyperplane field $\xi$, i.e. $\lambda|_{T{\mathcal L}}\equiv 0$.\\
Characteristic chords occur naturally in classical mechanics. Assume that $M$ is a compact hypersurface of contact type in the cotangent bundle of some manifold. Moreover, assume that $M$ is given as the zero set of a Hamiltonian function $H$ of the form $H(q,p)=\langle p,p\rangle+V(q)$, where $p$ represents coordinates along the fibers and $q$ stands for local coordinates on the base manifold, the pairing $\langle\,.\,,\,.\,\rangle$ denotes some bundle metric on the cotangent bundle and $V$ is a smooth function on the base manifold. If the intersection $\{p=0\}\cap M$ is a submanifold of $M$, then it is a Legendrian submanifold. Characteristic chords are trajectories of the given Hamiltonian system, where the momentum is zero at two different times, but the position $q$ of the system at these times is different. These are oscillating solutions, and they were investigated by Seifert in 1948 \cite{Seifert} and others since the 1970's \cite{Ambrosetti}, \cite{Bolotin}, \cite{VanGroesen}, \cite{Weinstein}.\\
In 1986, V.I. Arnold conjectured the existence of characteristic chords on the three dimensional sphere for any contact form inducing the standard contact structure and for any Legendrian knot \cite{Arnold}. After a partial result by the author in \cite{A3} this conjecture was finally confirmed by K. Mohnke in \cite{Mohnke}. It is natural to ask the existence question for characteristic chords not only for $M=S^3$, but also for general contact manifolds. A new invariant for Legendrian knots and contact manifolds proposed by Y. Eliashberg, A. Givental and H. Hofer in \cite{SFT} ('Relative Contact Homology') is actually based on counting characteristic chords and periodic orbits of the Reeb vector field.\\
In this paper we will use the method of filling by pseudoholomorphic curves which was developed in the previous papers \cite{part1}, \cite{part2} and \cite{part3} in order to prove a first existence result for characteristic chords on more general contact three manifolds. Before we can state the main result of this paper we need a few more definitions. Following \cite{Eliashberg89},\cite{Eliashberg92} we call $\xi$ 
overtwisted if there exists an embedded closed disk ${\cal D}$ in $M$
such that $T\partial \cal D\ \subset\ \xi \mid \cal D$ and $ \partial \cal D$
 does not contain any point $m$, where 
$T_m {\cal D}\ =\xi_{m}$. We call such a disk ${\cal D}$ an overtwisted disk. 
If such a disk does not exist we call the contact structure
tight. Similarly, a contact form is called overtwisted or tight if it induces a contact structure with the corresponding property. If ${\mathcal L}$ is a homologically trivial Legendrian knot in a three dimensional contact manifold then there are two 'classical' invariants, meaning invariants under Legendrian isotopy. There is the Thurston--Bennequin number $\mbox{tb}({\mathcal L})$ and the rotation number $r({\mathcal L},F)$, where $F$ is an embedded surface whose boundary is ${\mathcal L}$. Suppose ${\mathcal L}'$ is a knot obtained from ${\mathcal L}$ by slightly pushing ${\mathcal L}$ along some vector field transversal to $\xi$. The Thurston--Bennequin number $\mbox{tb}({\mathcal L})$ is then defined as the intersection number of ${\mathcal L}'$ with the surface $F$. The definition does not depend on the choice of the vector field and the spanning surface. For the definition of the rotation number, we have to assume that the Legendrian knot comes with a fixed orientation. Let then $v$ be a positively oriented, nowhere vanishing tangent vector field to ${\mathcal L}$. Let $\phi:\xi|_F\rightarrow F\times{\bf R}^2$ be a trivialization of the contact structure over the spanning surface $F$. The rotation number $r({\mathcal L},F)$ is then defined as the degree of the ${\bf R}^2$--component of $\phi\circ v$. This does not depend on the choice of the trivialization, but it depends on the relative homology class $[F]\in H_2(M,{\mathcal L})$.\\
We define a number $0<\inf(\lambda)\le\infty$ as follows:
\begin{equation}\label{definition-of-inf-number}
\inf(\lambda):=\inf\{T\,|\,x\ \mbox{is a T--periodic contractible orbit of}\ X_{\lambda}\}.
\end{equation}
If there are no contractible periodic orbits then $\inf(\lambda)=+\infty$. The number $\inf(\lambda)$ is obviously greater than zero since periodic orbits with very small period would lie in a Darboux--chart ($X_{\lambda}$ is bounded !) which is not possible. If $F$ is an oriented embedded surface which bounds a given Legendrian knot, then we define the '$\lambda$--volume' of $F$, $\mbox{vol}_{\lambda}(F)$ as follows: Let $\sigma$ be a 2--form on $F$ inducing the given orientation. We may define a function $g:F\rightarrow{\bf R}$ by the requirement $d\lambda|_F=g\cdot\sigma$. Then
\begin{equation}\label{definition-of-volume}
\mbox{vol}_{\lambda}(F):=\int_F|g|\,\sigma.
\end{equation}   
This definition does not depend on the choice of the form $\sigma$. We will prove the following existence result in this paper:
\begin{theorem}\label{main-result}
Let $M$ be a closed three dimensional manifold with a contact form $\lambda$ such that $\ker\lambda$ is a tight contact structure. Assume that ${\mathcal L}\subset M$ is a Legendrian knot which bounds an embedded disk ${\mathcal D}$ such that
\begin{equation}\label{lambda-volume-assumption}
\inf(\lambda)>\mbox{vol}_{\lambda}({\mathcal D})
\end{equation}
and $\mbox{tb}({\mathcal L})=-1$, $r({\mathcal L},{\mathcal D})=0$, i.e. the Legendrian has the largest possible Thurston--Bennequin number. Then there exists a characteristic chord for the Legendrian knot ${\mathcal L}$.
\end{theorem}

{\bf Remarks}:\\
The purpose of theorem \ref{main-result} is {\it not} to provide a tool for checking concrete examples. It is an existence result for characteristic chords which does not assume $M$ to be a particular three manifold like $S^3$ or ${\bf R}^3$ as it is the case with previous results. There is a large class of contact manifolds and Legendrian knots where theorem \ref{main-result} applies: Not every closed three dimensional contact manifold admits a tight contact structure \cite{Etnyre-Honda-1}, but there are several procedures known today for the construction of tight contact structures on a variety of manifolds (see for example \cite{Honda-1}, \cite{Honda-2}, \cite{Etnyre-Honda-2}, \cite{Etnyre-Ghrist} and the references cited in those papers). Moreover, in any closed tight contact three manifols there is an abundance of topological Legendrian unknots with the required Thurston--Bennequin and Rotation numbers (see \cite{Eliashberg-Fraser}). The assumption (\ref{lambda-volume-assumption}) is some kind of geometric 'smallness' assumption on the Legendrian knot.\\ 
The purpose of this paper is to give a first application of the filling method developed in \cite{part1}, \cite{part2} and \cite{part3} without any further technical complications. The assumption $\inf(\lambda)>\mbox{vol}_{\lambda}({\mathcal D})$ is vacuous if there are no contractible periodic orbits of the Reeb vector field. The purpose of this assumption is the following: The filling method used for proving theorem \ref{main-result} yields either a characteristic chord or a contractible periodic orbit with period bounded by $\mbox{vol}_{\lambda}({\mathcal D})$. The assumption (\ref{lambda-volume-assumption}) simply forbids the latter possibility. The author is convinced that there should be a more natural assumption in terms of Contact Homology to express the impact of periodic orbits on the Chord Problem.\\
In this paper we only consider Legendrian knots which are topological unknots and which have a certain Thurston--Bennequin number because we want to avoid dealing with saddle--type singular points of the characteristic foliation. Dealing with these issues significantly complicates matters because families of pseudoholomorphic strips cannot be extended further once they hit such a singular point. In fact, the boundary value problem discussed in this and the previous papers \cite{part1}, \cite{part2} and \cite{part3} would have to be modified in order to deal with these difficulties. These issues will be addressed in the forthcoming papers \cite{part5}, \cite{part6}.\\

The main tool of the proof are pseudoholomorphic curves in the symplectisation $({\bf R}\times M,d(e^t\lambda))$ of $M$. We are going to consider a special type of almost complex structures $\tilde{J}$ on ${\bf R}\times M$. We pick a complex structure $J:\xi\rightarrow\xi$ such that $d\lambda\circ(\mbox{Id}\times J)$ is a bundle metric on $\xi$. We then define an almost complex structure on ${\bf R}\times M$ by demanding $\tilde{J}\equiv J$ on $\xi$ and sending $\partial/\partial t$ (the generator of the ${\bf R}$--component) onto the Reeb vector field. Then $\tilde{J}(p)$ has to map $X_{\lambda}(p)$ onto $-\partial/\partial t$.\\
If $S$ is a Riemann surface with complex structure $j$ then we define a 
map 
$$
\tilde{u}=(a,u):S\longrightarrow{\bf R}\times M
$$
to be a pseudoholomorphic curve if
$$
D\tilde{u}(z)\circ j(z)=\tilde{J}(\tilde{u}(z))\circ D\tilde{u}(z)\ \mbox{for all}\ z\in S.
$$
If $(s,t)$ are conformal coordinates on $S$ then this becomes:
\[
\partial_s\tilde{u}+\tilde{J}(\tilde{u})\partial_t\tilde{u}=0.
\]
We are interested only in pseudoholomorphic curves which have finite energy in the sense that
$$
E(\tilde{u}):=\sup_{\phi\in\Sigma}\int_{S}
\tilde{u}^{\ast}d(\phi\lambda) <+\infty
$$
where $\Sigma:=\{\phi\in C^{\infty}({\bf R},[0,1])\,|\,\phi'\ge 
0\}$. The following results show that nontrivial pseudoholomorphic planes and half--planes with finite energy lead to the existence of periodic orbits and characteristic chords respectively \cite{Hofer-Weinstein-conj}:
\begin{theorem}\label{finite energy plane}
Assume $M$ is a closed manifold with contact form $\lambda$, 
and let $\tilde{u}=(a,u)$ be a non constant pseudoholomorphic plane with finite energy. Then
$$
T:=\int_{{\bf C}}u^{\ast}d\lambda>0
$$
and every sequence $(R_k')_{k\in{\bf N}}$ of positive real numbers 
tending to infinity has a subsequence 
$(R_k)_{k\in{\bf N}}$ so that ~$ u(R_ke^{\frac{2\pi i}{T}\cdot})$ converges 
in $C^{\infty}({\bf R}/T{\bf Z},M)$ to some $x$ which satisfies
$$
\dot{x}(t)=X_{\lambda}(x(t))\ \mbox{ and }\ \ x(0)=x(T).
$$
\end{theorem}\qed \\
There is a corresponding result for characteristic chords. If $H^+:=\{s+it\in{\bf C}\,|\,t\ge 0\}$ is the closed upper half of the 
complex plane, $M$ a manifold with contact form $\lambda$, and 
${\cal L}\subset M$ a Legendrian submanifold (or a Legendrian knot if $M$ is three--dimensional) then we define a finite energy half--plane to be a map
$$
\tilde{u}=(a,u):H^+\longrightarrow{\bf R}\times M
$$
that satisfies the following conditions:
\begin{enumerate}
	\item  $\partial_s\tilde{u}+\tilde{J}(\tilde{u})\partial_t\tilde{u}=0$ 
	on $\overset{\circ\ }{H^+}$

	\item  $\tilde{u}(\partial H^+)\subset {\bf R}\times{\cal L}$

	\item  $u(H^+)$ is contained in a compact region $K\subset M$

	\item  $\tilde{u}$ has finite energy: 
	$E(\tilde{u})<+\infty.$
\end{enumerate}
Existence of a nontrivial finite energy half--plane implies existence of a characteristic chord (see \cite{A1} for a proof):
\begin{theorem}\label{finite energy half--plane}
Let $\tilde{u}$ be a finite energy half--plane which is in 
addition non constant. Then 
$T:=\int_{H^+}u^{\ast}d\lambda>0$ and any sequence of positive 
real numbers tending to $+\infty$ has a subsequence 
$R_k\rightarrow+\infty$,
 so that the maps 
$$[0,T]\longrightarrow M$$
$$t\longmapsto u(R_k\,e^{\pi i\frac{t}{T}})$$
converge in $C^{\infty}$ to some orbit $x$ of the Reeb vector field 
$X_{\lambda}$ with $x(0) , x(T) \in{\cal L}$.
\end{theorem}\qed \\
The aim of this paper is to show the existence of a nontrivial finite energy half--plane under the assumptions of theorem \ref{main-result}. If $F$ is an embedded oriented surface in a three dimensional contact manifold $(M,\xi=\ker\lambda)$ then there is a distinguished singular foliation on $F$, the so--called characteristic foliation. The singularities consist of the points where the contact structure is tangent to the surface. The leaves of the foliation are the curves on $F$ which are tangent to $\ker\lambda$. A singularity $p\in F$ is called positive if the orientation of $(T_pF,d\lambda_p)$ coincides with the prescribed orientation of $F$, and negative otherwise. If $\sigma$ is a 2--form on $F$ inducing the given orientation then the vector field $Z$ defined by $i_Z\sigma=j^{\ast}\lambda$ induces the characteristic foliation on $F$ ($j:F\hookrightarrow M$ is the inclusion). The zeros of $Z$ then correspond to the singularities of the characteristic foliation. If $Z(p)=0$ is a saddle--type singularity then we call $p$ a hyperbolic singularity. In the case of a sink or a source we call $p$ an elliptic singularity.
\section{The filling method}
The strategy of the proof is the following: In order to show existence of a characteristic chord it is sufficient to show the existence of a nontrivial finite energy half--plane as in theorem \ref{finite energy half--plane}. We will accomplish this by considering a different boundary value problem involving the nonlinear Cauchy Riemann equation. This boundary value problem is set up in such a way that the moduli space of its solutions is not compact. A bubbling off analysis then yields existence of either a finite energy plane or a finite energy half--plane. Assumption (\ref{lambda-volume-assumption}) of theorem \ref{main-result} excludes the case of a finite energy plane proving the existence of a characteristic chord.\\
Assume in this section that $M$ is a closed three dimensional manifold with contact form $\lambda$ and that ${\mathcal L}$ is a homologically trivial Legendrian knot with Seifert surface ${\mathcal D}$. Moreover, let $e\in{\mathcal L}$ be an elliptic singular point. We consider the following boundary value problem on the infinite strip $S={\bf R}\times [0,1]$:
\begin{equation}\label{main-boundary-value-problem}
 \begin{array}{ll}
          \tilde{u}=(a,u):S\longrightarrow\RM & \\
          \partial_s\tu+\tilde{J}(\tu)\partial_t\tu=0 & \\
          \tu(s,0)\subset{\bf R}\times{\cal L} & \\
          \tu(s,1)\subset\{0\}\times{\cal D}^{\ast} & \\
          u(0,0)=e & \\
           0< E(\tu)<+\infty & 
          \end{array},
\end{equation}
where ${\mathcal D}^{\ast}$ denotes the Seifert surface without its singular points. Pseudoholomorphic strips with mixed (Lagrangian) boundary conditions were also studied by A. Floer in his theory of Lagrangian intersections \cite{Floer} but there are important differences. In contrast to A. Floer's theory, the boundary conditions $\{0\}\times{\mathcal D}^{\ast}$ and ${\bf R}\times{\mathcal L}$ have non--transverse intersection, and $\{0\}\times{\mathcal D}^{\ast}$ is not Lagrangian. Moreover, the 'ends' $\lim_{s\rightarrow\pm\infty}\tu(s,t)$ of the solutions of the problem (\ref{main-boundary-value-problem}) are not fixed, they are allowed to slide along ${\mathcal L}$. In addition, theorem \ref{local-existence-theorem} below requires a particular almost complex structure on a neighborhood of the elliptic singular points on the boundary. A. Floer's proof of transversality by choosing a generic almost complex structure does not work here. With a bit more effort (using also that ${\bf R}\times M$ is four dimensional) one can show that transversality holds automatically (see \cite{part2}). The previous papers \cite{part1}, \cite{part2} and \cite{part3} take care of all these analytic issues. The purpose of this section is to give a brief summary. 

\begin{theorem}\label{local-existence-theorem}
{\bf (Local existence of solutions)}\\Let $(M,\lambda)$ be a three dimensional contact manifold. Moreover, let ${\mathcal L}$ be a Legendrian knot which bounds an embedded surface ${\mathcal D}'$ so that the characteristic foliation has only finitely many singular points. Then there is another embedded surface ${\mathcal D}$ which is a smooth $C^0$--small perturbation of ${\mathcal D}'$ having the same boundary and the same singular points as ${\mathcal D}'$ and a $d\lambda$--compatible complex structure $J:\ker\lambda\rightarrow\ker\lambda$ so that the following is true: Near each elliptic singular point $e\in\partial{\mathcal D}={\mathcal L}$ there are embedded solutions $\tu_{\tau}$ , $0<\tau<1$ to the boundary value problem (\ref{main-boundary-value-problem}) with the properties:
\begin{itemize}
\item $\tu_{\tau}(S)\cap\tu_{\tau'}(S)=\emptyset$ if $\tau\neq\tau'$,
\item $\tu_{\tau}\rightarrow e$ uniformly with all derivatives as $\tau\rightarrow 0$,
\item the family $\tu_{\tau}$ depends smoothly on the parameter $\tau$,
\item each map $u_{\tau}$ is transverse to the Reeb vector field, i.e. $\pil\pas u_{\tau}(z)\neq 0$ for all $z\in S$,
\item The Maslov indices $\mu(\tu_{\tau})$ all equal zero.
\end{itemize}
\end{theorem}

{\bf Proof:} See \cite{part3}.
\qed\\

The Maslov--index $\mu(\tu_{\tau})$ is a topological quantity associated to the boundary condition $\{0\}\times{\mathcal D}$. Its precise definition is not relevant for this paper, see \cite{part2}.

\begin{proposition}\label{7.1.}
Let $\tu=(a,u):S\longrightarrow\RM$ be a non--constant solution of
the boundary value problem (\ref{main-boundary-value-problem}). 
\begin{itemize}
\item Then the path
$s\longmapsto u(s,1)$ is transverse to the characteristic
foliation, i.e. $\pas u(s,1)\not\in\ker\lambda(u(s,1))$. We actually have
\[
0<\pat a(s,1)=-\lambda(u(s,1))\pas u(s,1)
\]
for all $s\in{\bf R}$.
\item We have $a(s,t)<0$ whenever $0\le t<1$,
\item The pseudoholomorphic strip never hits $\{0\}\times{\mathcal L}$, i.e.
\[
\tu(S)\cap (\{0\}\times{\mathcal L})=\emptyset.
\]
In particular,
\[
\lim_{s\rightarrow\pm\infty}\tu(s,t)\not\in\tu(S).
\]
\end{itemize}
\end{proposition}

{\bf Proof:} See \cite{part2}.
\qed\\

The subject of the paper \cite{part1} is the asymptotic behavior of solutions to (\ref{main-boundary-value-problem}) for $s\rightarrow\pm\infty$. Among other things we have shown that for any solution $\tu$ of (\ref{main-boundary-value-problem}) the limits $p_{\pm}=\lim_{s\rightarrow\pm\infty}\tu(s,t)\in\{0\}\times {\mathcal L}$ exist. Moreover, in local coordinates the function $\tu(s,t)-p_{\pm}$ and all its derivatives converge to zero like $e^{\lambda_{\pm}s}$, where $\lambda_+<0$, $\lambda_->0$ are integer multiples of $\pi/2$ (see \cite{part1}). The local solutions in theorem \ref{local-existence-theorem} actually decay at the rate $e^{-\pi/2|s|}$ (see \cite{part3}). The main result of the paper \cite{part2} states that existence of a suitable embedded solution implies the existence of a whole 1--parameter family of embedded solutions near by:

\begin{theorem}\label{main-implicit-function-theorem}
{\bf (Implicit function theorem)}\\Let $\tu_0=(a_0,u_0)$ be an embedded solution of (\ref{main-boundary-value-problem}) so that its Maslov--index $\mu(\tu_0)$ vanishes. Assume moreover, that $|\tu_0(s,t)-p_{\pm}|$ decays either like $e^{-\pi|s|}\,$\footnote{Solutions which decay like $e^{-\pi|s|}$ or faster actually do not occur, as the Compactness Result below shows. We nevertheless have to deal with solutions of decay $e^{-\pi|s|}$ in the Implicit Function Theorem because the proof of the Compactness Result in \cite{part3} requires us to do so in order to show their non--existence.}  or like $e^{-\frac{\pi}{2}|s|}$ for large $|s|$ in local coordinates near the points $p_{\pm}:=\lim_{s\rightarrow\pm\infty}\tu_0(s,t)$ and that $p_-\neq p_+$. Assume also that
\[
\mbox{dist}(u_0({\bf R}\times\{1\}),\Gamma)>0,
\]
where $\Gamma$ denotes the set of the singular points on the Seifert surface ${\mathcal D}$. Then there is a smooth family $(\tilde{v}_{\tau})_{-1<\tau<1}$ of embedded solutions of (\ref{main-boundary-value-problem}) with the following properties:
\begin{itemize}
\item $\tilde{v}_0=\tu_0$,
\item The solutions $\tilde{v}_{\tau}$ have the same Maslov--index and the same decay rates as $\tu_0$,
\item The sets 
\[
U_{\pm}:=\bigcup_{-1<\tau< 1}\{\lim_{s\rightarrow\pm\infty}\tilde{v}_{\tau}(s,t)\}
\]
are open neighborhoods of the points $p_{\pm}$ in ${\mathcal L}$.
\end{itemize}
If $|\tu_0(s,t)-p_{\pm}|$ decays like $e^{-\frac{\pi}{2}|s|}$ for both $s\rightarrow+\infty$ and $s\rightarrow-\infty$ then we have in addition
\begin{itemize}
\item $\tilde{v}_{\tau}(S)\cap\tilde{v}_{\tau'}(S)=\emptyset$ if $\tau\neq\tau'$.
\end{itemize}
\end{theorem}

{\bf Proof:} See \cite{part2}.
\qed\\

The assumption on the decay of $\tu_0$ and its Maslov--index are important since the Fredholm index depends on this data. This means that we need an enhanced compactness result which guarantees the correct data plus the embedding property for a $C^{\infty}_{loc}$--limit of a sequence of solutions. Analytic arguments yield the exponential decay of the limit while the specific decay rate and the embedding property follows from topological arguments. We have shown the following compactness result in the paper \cite{part3}:

\begin{theorem}\label{compactness-result}
{\bf (Compactness result)}\\Let $(\tu_{\tau})_{0\le\tau<\tau_0}=(a_{\tau},u_{\tau})_{0\le\tau<\tau_0}$ be a smooth family of embedded solutions to the boundary value problem
\[
\begin{array}{ll}
          \tilde{u}=(a,u):S\longrightarrow\RM & \\
          \partial_s\tu+\tilde{J}(\tu)\partial_t\tu=0 & \\
          \tu(s,0)\subset{\bf R}\times{\cal L} & \\
          \tu(s,1)\subset\{0\}\times{\cal D}\backslash\Gamma & \\
          u(0,0)=e & \\
           0< E(\tu)<+\infty & 
          \end{array},
\]
where ${\mathcal D}\subset M$ is an embedded surface bounding the Legendrian knot ${\mathcal L}$, $\Gamma\subset {\mathcal D}$ is the set of singular points and $e\in\Gamma\cap{\mathcal L}$ is an elliptic singular point on the boundary of ${\mathcal D}$. We impose the following conditions on the solutions $\tu_{\tau}$:
\begin{itemize}
\item $\tu_{\tau'}(S)\cap\tu_{\tau''}(S)=\emptyset$ if $\tau'\neq \tau''$,
\item 
\[
\mbox{dist}\left(\bigcup_{0<\delta\le\tau<\tau_0}\{u_{\tau}(s,1)\,|\,s\in{\bf R}\}\,,\,\Gamma\right)>0,
\]
\item For small $\tau$ the solutions $\tu_{\tau}$ coincide with the local solutions of theorem \ref{local-existence-theorem} near $e$, in particular, they decay like $e^{-\pi/2|s|}$,
\item There is a uniform gradient bound, i.e.
\[
\sup_{0\le\tau<\tau_0}\|\nabla\tu_{\tau}\|_{C^0(S)}<\infty.
\]
\end{itemize}
Then for every sequence $\tau'_k\nearrow\tau_0$ there is a subsequence $\tau_k$ such that the family $\tu_{\tau_k}$ converges in $C^{\infty}_{loc}$ to another solution (as $k\rightarrow\infty$) $\tu_{\tau_0}$ with finite energy such that also $\mbox{dist}\big(\{u_{\tau_0}(s,1)\,|\,s\in{\bf R}\}\,,\,\Gamma\big)>0$. Moreover,
\begin{enumerate}
\item every sequence $\tau_k$ yields the same limit, i.e. $\tu_{\tau_0}=\lim_{\tau\nearrow\tau_0}\tu_{\tau}$, and the convergence is uniform on $S$ with all derivatives,
\item $\tu_{\tau_0}$ is an embedding,
\item the Maslov--index $\mu(\tu_{\tau_0})$ of $\tu_{\tau_0}$ equals $0$,
\item The solution $\tu_{\tau_0}(s,t)$ has the same rate of decay for large $|s|$ as the maps $\tu_{\tau}$, i.e. $|\lambda_{\pm}|=\frac{\pi}{2}$,
\item $\tu_{\tau_0}(S)\cap\tu_{\tau}(S)=\emptyset$ for all $0\le\tau<\tau_0$.
\end{enumerate}
\end{theorem}

{\bf Proof:} See \cite{part3}.
\qed\\

Assume $\tu_{\tau}=(a_{\tau},u_{\tau})$ , $-1<\tau\le 0$ is a continuous family of embedded pseudoholomorphic strips as in (\ref{main-boundary-value-problem}) with pairwise disjoint images. Let $\tilde{v}=(b,v)$ be another embedded solution of the boundary value problem (\ref{main-boundary-value-problem}).
In the paper \cite{part3} we have studied the intersection properties of the pseudoholomorphic curve $\tilde{v}$ with the family $\tu_{\tau}$. The following theorem states that 'there is no isolated first intersection'.

\begin{theorem}
\label{intersection-of-families-0}
{\bf (No isolated first intersection)}\\Assume $\tu_{\tau}=(a_{\tau},u_{\tau})$ , $-1<\tau\le 0$ is a smooth family of embedded solutions of (\ref{main-boundary-value-problem}) with pairwise disjoint images and let $\tilde{v}=(b,v)$ be another embedded solution. Moreover, we assume that $\tilde{u}_{\tau}$ and $\tilde{v}$ have disjoint images for $\tau<0$, but the intersection of $\tu_0({\bf R}\times [0,1])$ with $\tilde{v}({\bf R}\times [0,1])$ is not empty. Then the image of $\tilde{v}$ is contained in the image of $\tu_0$ or vice versa except in the following case: If the first intersection occurs at the boundary ${\bf R}\times\{0\}$, i.e. if $\tu_0(p)=\tilde{v}(q)$ for $p,q\in{\bf R}\times\{0\}$, and $\pas u_0(p)$, $\pas v(q)\in T_{u_0(p)}{\mathcal L}$ do not have the same orientation then we can only conclude that $\tu_0({\bf R}\times \{0\})=\tilde{v}({\bf R}\times\{0\})$.
\end{theorem}
{\bf Proof:} See \cite{part3}.
\qed\\

{\bf Remark:} In most cases it is sufficient to know that the images of $\tu_0$ and $\tilde{v}$ agree along the boundary ${\bf R}\times\{0\}$.\\

We also have to exclude the situation where the intersection occurs at infinity.
\begin{theorem}\label{intersections-at-infinity}
Let $\tu_{\tau}=(a_{\tau},u_{\tau})$ , $-1<\tau\le 0$ be a smooth family of embedded solutions of (\ref{main-boundary-value-problem}) with pairwise disjoint images. Let $\tilde{v}$ be another embedded solution. Assume that all the maps $\tu_{\tau}$ and $\tilde{v}$ have the same exponential decay rate $\lambda_+=-\pi/2$ as $s\rightarrow+\infty$. We assume also that $\tilde{v}$ and $\tu_0$ converge to the same point on $\{0\}\times{\mathcal L}$ as $s\rightarrow+\infty$, but 
\[
\tilde{u}_{\tau}(S)\cap \tilde{v}(S)=\emptyset\ \mbox{for all}\ \tau<0
\]
Then 
\[
\tilde{v}({\bf R}\times\{0\})=\tu_{0}({\bf R}\times\{0\}).
\] 
\end{theorem}
{\bf Proof:} See \cite{part3}. \qed

\begin{figure}[!ht]
\begin{center}
\epsfig{file=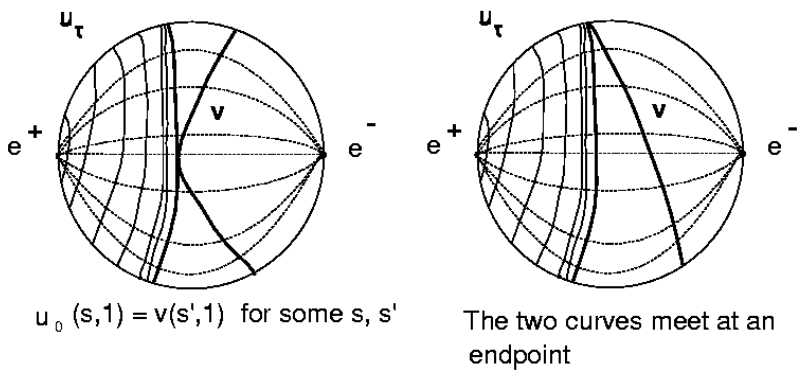,height=5cm,angle=0}
\caption{Intersection scenarios forbidden by theorems \ref{intersection-of-families-0} and \ref{intersections-at-infinity}.}\label{intersections}
\end{center}
\end{figure}

In figure \ref{intersections} we visualize the solutions to the boundary value problem (\ref{main-boundary-value-problem}) by looking at $u_{\tau}({\bf R}\times\{1\})$ and $\tilde{v}({\bf R}\times\{1\})$ which are all curves on the Seifert surface ${\mathcal D}$.

\section{Manipulating the characteristic foliation}

\begin{figure}[!ht]
\begin{center}
\epsfig{file=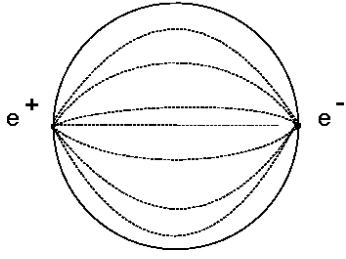,height=4cm,angle=0}
\caption{Characteristic foliation in the case where $\mbox{tb}({\mathcal L})=-1$ and $r({\mathcal L},{\mathcal D})=0$.}\label{special-elliptic-form}
\end{center}
\end{figure}

The results of the previous section were all formulated without any particular assumptions on the Legendrian knot ${\mathcal L}$ apart from being homologically trivial. We will now return to the situation of theorem \ref{main-result} where the Seifert surface ${\mathcal D}$ is a disk, the Legendrian has particular Thurston-Bennequin and rotation numbers and $M$ is a tight contact manifold. Because ${\mathcal L}=\partial{\mathcal D}$ is Legendrian, the boundary of the disk belongs to the characteristic foliation and because $\ker\lambda$ is tight there must be singular points on the boundary. Following the paper \cite{Eliashberg-Fraser}, we may deform the disk ${\mathcal D}$ to another disk ${\mathcal D}'$ which satisfies the following conditions
\begin{itemize}
\item $\partial{\mathcal D}=\partial{\mathcal D}'={\mathcal L}$,
\item $|\mbox{vol}_{\lambda}({\mathcal D})-\mbox{vol}_{\lambda}({\mathcal D}')|$ is as small as we wish
\item the characteristic foliation on the new disk ${\mathcal D}'$ has exactly two singular points. They are both elliptic, and they both lie on the boundary (see figure \ref{special-elliptic-form}).
\end{itemize}

\section{Proof of theorem \ref{main-result}}

Near each boundary elliptic point $e$ there are solutions $\tu_{\tau}=(a_{\tau},u_{\tau})$ of (\ref{main-boundary-value-problem}) starting out with the family given by the local existence theorem, theorem \ref{local-existence-theorem}. By our intersection results, theorems \ref{intersection-of-families-0} and \ref{intersections-at-infinity}, the two solution families $\tu^+_{\tau}$, $\tu^-_{\tau}$ originating from $e^+$ and $e^-$ respectively can not intersect. Otherwise, the restrictions of the two families to ${\bf R}\times\{0\}$ would blend together into one family. If the families intersected then we could construct a smooth map $\Psi:Q\rightarrow {\bf R}\times {\mathcal L}$ defined on some closed rectangle $Q\subset{\bf R}^2$ so that its restriction to the boundary has degree one. We define a subset ${\mathcal O}\subset{\mathcal D}$, the set of obstacles, to be the union of the curves $u^-_{\tau}({\bf R}\times\{1\})$ including their end points $\lim_{s\rightarrow\pm\infty}u_{\tau}^-(s,t)\in{\mathcal L}$. The set of obstacles consists of objects on the Seifert surface which can not be hit by the family of curves $u^+_{\tau}({\bf R}\times\{1\})$. \\
We claim that there is a family of embedded solutions $\tu^+_{\tau}$, $0<\tau<1$ of (\ref{main-boundary-value-problem}) starting out from the elliptic singular point $e^+$ on the boundary as in theorem \ref{local-existence-theorem} such that 
\begin{equation}\label{grad-blow-up}
\sup_{0\le\tau<1}\|\nabla\tu^+_{\tau}\|_{C^0(S)}=+\infty.
\end{equation}
Let us prove this indirectly. Let ${\mathcal F}_e$ be the set of all smooth families of embedded solutions $(\tu_{\tau'})_{0<\tau'<\tau'_0}$ to (\ref{main-boundary-value-problem}) with $u_0=e^+\in{\mathcal L}$, which agree with the local family given by theorem \ref{local-existence-theorem} if $\tau'$ is sufficiently small.
Arguing indirectly, we assume that for each family $(\tu_{\tau'})\in{\mathcal F}_e$: 
\[
\sup_{0<\tau'<\tau'_0}\|\nabla\tu_{\tau'}\|_{C^0(S)}<+\infty.
\]
Applying the compactness result, theorem \ref{compactness-result}, any family $(\tu_{\tau'})_{0<\tau'<\tau'_0}\in{\mathcal F}_e$ extends to a family $(\tu_{\tau'})_{0<\tau'\le\tau'_0}\in{\mathcal F}_e$. Then $\tu_{\tau'_0}$ is also embedded, and it has the same decay rate and Maslov index as the rest of the family, namely $|\lambda_{\pm}|=\frac{\pi}{2}$. Moreover, the image of $\tu_{\tau'_0}$ is disjoint from the images of the $\tu_{\tau'}$ for $\tau'<\tau'_0$.\\  
We can now apply the implicit function  theorem, theorem \ref{main-implicit-function-theorem}, to the solution $\tu_{\tau_0}$. We obtain a larger family $(\tu_{\tau'})_{0<\tau'<\tau'_0+\ve}\in{\mathcal F}_e$, and the boundary curves on ${\mathcal D}$ move closer to the obstacles ${\mathcal O}$:
\[
\bigcup_{0<\tau'<\tau'_0}u_{\tau'}({\bf R}\times\{1\}) \stackrel{\neq}{\subset}\bigcup_{0<\tau'<\tau'_0+\ve}u_{\tau'}({\bf R}\times \{1\}).
\]
By assumption, we have again a gradient bound for this larger family
\[
\sup_{0<\tau'<\tau'_0+\ve}\|\nabla\tu_{\tau'}\|_{C^0(S)}<+\infty.
\]
On the other hand there is a positive constant $\delta$ such that 
\[
\mbox{dist}\left(\bigcup_{\tau}u_{\tau}({\bf R}\times\{1\})\,,\,{\mathcal O}\right)\ge\delta
\]
for any family $(\tu_{\tau})\in{\mathcal F}_e$. If not, then we could find a family of solutions in ${\mathcal F}_e$ which intersects the set of obstacles ${\mathcal O}$. But we have just seen that this is impossible. This proves our claim (\ref{grad-blow-up}).\\
If $\tu_k$ is a sequence of solutions of (\ref{main-boundary-value-problem}) and $z_k\in S$ with $R_k:=|\nabla\tu_k(z_k)|\rightarrow \infty$ then we use a standard rescaling argument (for example as in section 3 of the paper \cite{part1}). Depending on where bubbling occurs we either obtain a pseudoholomorphic plane, a half--plane with boundary condition ${\bf R}\times{\mathcal L}$ or a half--plane with boundary condition $\{0\}\times{\mathcal D}^{\ast}$ where the boundary curve has positive distance from the set ${\mathcal O}$ and is everywhere transverse to the characteristic foliation which is also impossible. In all cases the resulting pseudoholomorphic curve is not constant and has energy bounded by $\mbox{vol}_{\lambda}({\mathcal D})$ because the energy of any solution $\tu_{\tau}$ of (\ref{main-boundary-value-problem}) is bounded by $\mbox{vol}_{\lambda}({\mathcal D})$ (see proposition 2.3 in \cite{part2}). In the case of a pseudoholomorphic plane we would obtain a contractible periodic orbit of the Reeb vector field (theorem \ref{finite energy plane}) with period no larger than $\mbox{vol}_{\lambda}({\mathcal D})$, which is impossible because of the assumption 
\[
\inf(\lambda)>\mbox{vol}_{\lambda}({\mathcal D}).
\]
The case of a half--plane with boundary condition $\{0\}\times {\mathcal D}$ is also impossible since we could remove the singularity at infinity and obtain a disk with boundary curve on ${\mathcal D}^{\ast}$ and transverse to the characteristic foliation. The only remaining possibility then yields a characteristic chord. In view of the assumption
\[
\inf(\lambda)>\mbox{vol}_{\lambda}({\mathcal D})
\]
we obtain a genuine chord $x(t)$, i.e. $x(0),x(T)\in{\mathcal L}$ for some $T>0$ and $x(0)\neq x(T)$. This completes the proof of theorem \ref{main-result}.
\qed\\

\end{document}